\documentclass[a4paper,UKenglish]{lipics_md}
\usepackage{microtype}
\usepackage{mathtools}
\usepackage[normalem]{ulem}

\title{Poisson--Delaunay Mosaics of Order $k$\footnote{This work is partially
     supported by the Austrian Science Fund (FWF), grant no.\ I02979-N35.}}
\titlerunning{Poisson--Delaunay Mosaics of Order \textit{k}}

\author[]{Herbert Edelsbrunner}
\author[]{Anton Nikitenko}
\affil[]{IST Austria (Institute of Science and Technology Austria),
  Am Campus 1, \\ 3400 Klosterneuburg, Austria,
  \texttt{edels@ist.ac.at}, \texttt{anton.nikitenko@ist.ac.at}}
\authorrunning{H. Edelsbrunner and A. Nikitenko}

\mathsubjclass{60D05 Geometric probability and stochastic geometry,
  68U05 Computer graphics; computational geometry.} 

\newcommand{\mm}[1] {\ifmmode{#1}\else{\mbox{\(#1\)}}\fi}

\newcommand{\ignore}[1]{}

\newcommand{\ourproof}{\begin{proof}}
\newcommand{\eop}{\end{proof}}  %

\newcommand{\Rspace}       {\mm{{\mathbb R}}}
\newcommand{\EE}[1]        {\mm{{\mathbb E}{[{#1}]}}} 
\newcommand{\PP}[2]        {\mm{{\mathbb P}_{#1}{[{#2}]}}}
\newcommand{\One}[1]       {\mm{{{\bf 1}_{#1}}}}

\newcommand{\Gama}[1]      {\mm{{\Gamma}{\left({#1}\right)}}}
\newcommand{\iGama}[2]     {\mm{{\gamma}{\left({#1};\,{#2}\right)}}}
\newcommand{\Del}[2]       {\mm{{\rm Del}_{#1}{({#2})}}}
\newcommand{\Vor}[2]       {\mm{{\rm Vor}_{#1}{({#2})}}}
\newcommand{\domain}[1]    {\mm{{\rm dom}{({#1})}}}
\newcommand{\power}[1]     {\mm{{\pi}_{#1}}}
\newcommand{\Uin}          {\mm{U_{\it \!in}}}
\newcommand{\Uon}          {\mm{U_{\it \!on}}}
\newcommand{\Uout}         {\mm{U_{\it \!out}}}
\newcommand{\Rfun}         {\mm{{\cal R}}}
\newcommand{\density}      {\mm{\rho}}
\newcommand{\xxx}          {\mm{{\bf x}}}

\newcommand{\Area}[3]      {\mm{\eta_{{#1}}^{{#2},{#3}}}}
\newcommand{\NFaces}[4]    {\mm{N}_{{#1},{#2}}^{{#3}}{({#4})}}
\newcommand{\cNumNB}[5]    {\mm{s_{{#1},{#2},{#3}}^{{#4},{#5}}}}
\newcommand{\cNum}[5]      {\mm{c_{{#1},{#2},{#3}}^{{#4},{#5}}}}
\newcommand{\dNum}[3]      {\mm{d_{{#1}}^{{#2},{#3}}}}
\newcommand{\CConst}[3]    {\mm{C_{{#1},{#2}}^{#3}}}

\newcommand{\Sphere}[2]    {\mm{\Sigma_{#1}{({#2})}}}
\newcommand{\radius}[2]    {\mm{r_{#1}{({#2})}}}
\newcommand{\SphereOnly}   {\mm{\Sigma}}
\newcommand{\Inn}[2]       {\mm{{\rm In}{({#2})}}}
\newcommand{\Onn}[2]       {\mm{{\rm On}{({#2})}}}
\newcommand{\inn}[2]       {\mm{{\rm in}{({#2})}}}
\newcommand{\onn}[2]       {\mm{{\rm on}{({#2})}}}

\newcommand{\conv}[1]      {\mm{\rm conv\,}{#1}}
\newcommand{\interior}[1]  {\mm{\rm int\,}{#1}}

\newcommand{\diff}         {\mm{\rm \,d}}
\newcommand{\norm}[1]      {\mm{\|{#1}\|}}
\newcommand{\card}[1]      {|{#1}|}
\newcommand{\Edist}[2]     {\mm{\|{#1}-{#2}\|}}

\newcommand{\nnn}          {\mm{\bf n}}

\newcommand{\Region}       {\mm{\Omega}}

\newcommand{\DimVor}       {\mm{\ell}}
\newcommand{\DimDel}       {\mm{j}}

\newcommand{\ourparagraph}[1] {\vspace{0.1in} \noindent \textbf{#1}}
\newcommand{\Skip}[1]      {}

\begin{document}
\maketitle

\begin{abstract}
  The order-$k$ Voronoi tessellation of a locally finite set
  $X \subseteq \Rspace^n$ decomposes $\Rspace^n$ into convex domains
  whose points have the same $k$ nearest neighbors in $X$.
  Assuming $X$ is a stationary Poisson point process,
  we give explicit formulas for the expected number
  and total area of faces of a given dimension per unit volume of space.
  We also develop a relaxed version of discrete Morse theory and generalize
  by counting only faces, for which the $k$ nearest points
  in $X$ are within a given distance threshold.
\end{abstract}

\section{Introduction}
\label{sec:1}

\begin{figure}[hbt]
  \centering \resizebox{!}{1.8in}{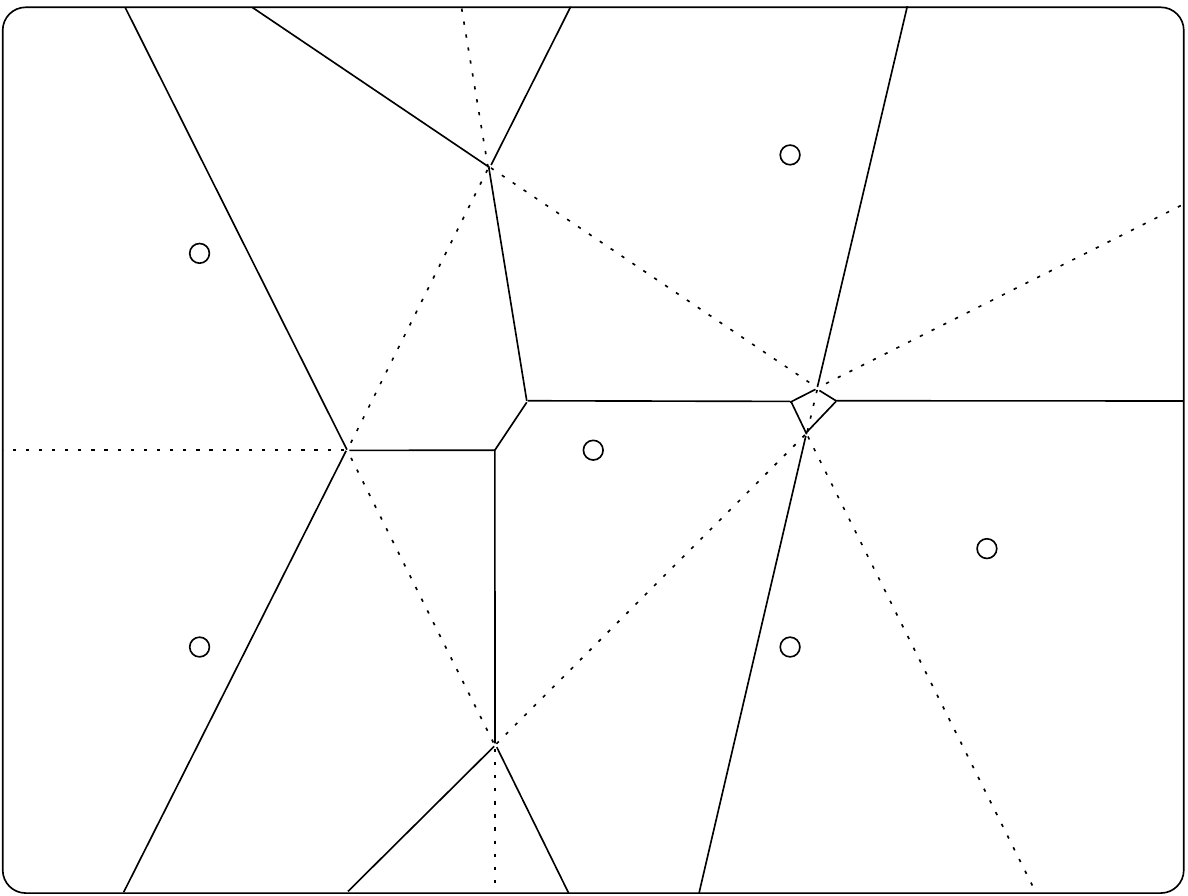}
  \caption{The dotted edges decompose the plane into the order-$1$ Voronoi domains,
    while the solid edges decompose it into order-$2$ Voronoi domains.
    Observe that the two tessellations share some of their vertices but not all.}
  \label{fig:Voronoi}
\end{figure}

Let $X\subseteq \Rspace^n$ be locally finite.
The \emph{Voronoi domain} of a subset $Q \subseteq X$,
denoted $\domain{Q}$, is the set of points $p \in \Rspace^n$
for which $\Edist{p}{x} \leq \Edist{p}{y}$
for all $x \in Q$ and all $y \in X \setminus Q$.
The \emph{order} of the domain is the cardinality of $Q$.
For any integer $k \geq 1$, the \emph{order-$k$ Voronoi tessellation}
of $X$ is the collection of order-$k$ Voronoi domains;
that is: domains of sets $Q \subseteq X$ with $\card{Q} = k$;
see \cite{Fej76,Lee82,ShHo75}.
Figure \ref{fig:Voronoi} illustrates this concept by superimposing
two tessellations of a finite set in the plane.
For $k = 1$, we get what is usually called the Voronoi diagram
or Voronoi tessellation \cite{Aur91}, which is
generically \emph{primitive} (or \emph{normal}).
This means that the common intersection of $\DimVor+1$ domains is either
empty or has dimension $n-\DimVor$.
As we will explain in Section \ref{sec:3},
this is also true for Voronoi tessellations of order $k$ in dimension $2$
but not for Voronoi tessellations of order $k \geq 2$ in dimension $n \geq 3$.

\begin{figure}[hbt]
  \centering \resizebox{!}{1.8in}{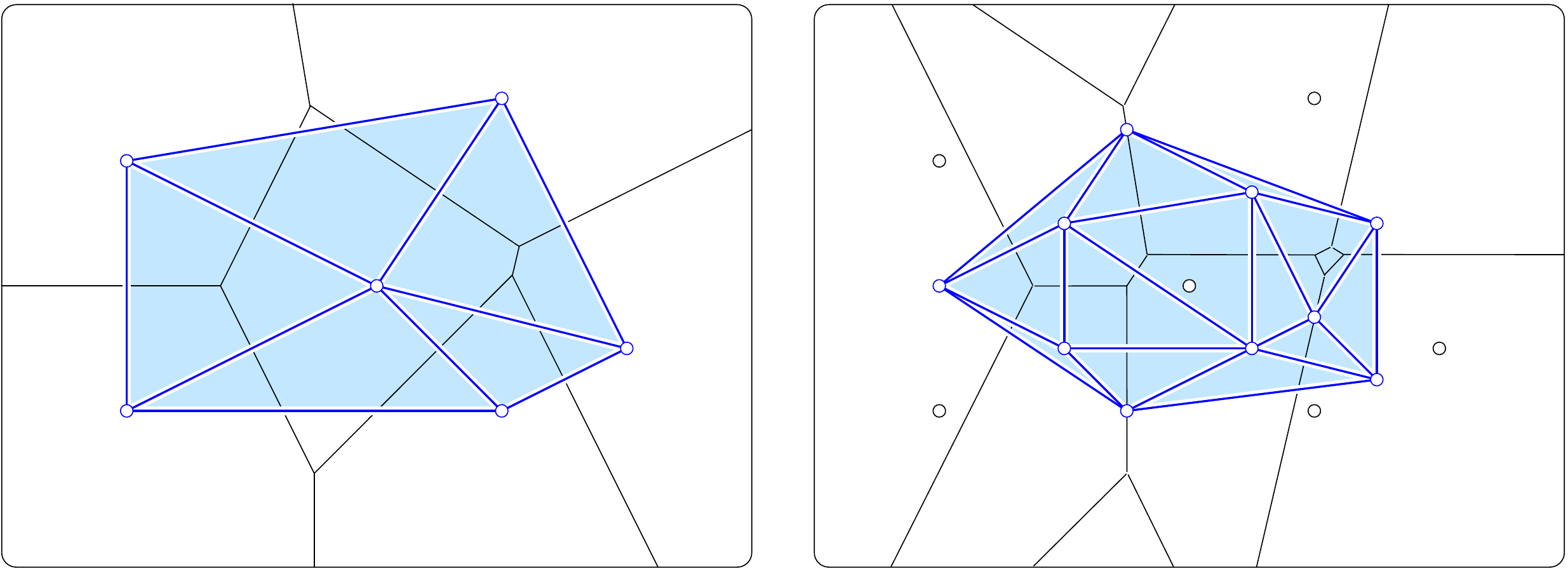}
  \caption{The order-$1$ Delaunay mosaic on the \emph{left}
    and the order-$2$ Delaunay mosaic on the \emph{right},
    both superimposed on their corresponding Voronoi tessellations.}
  \label{fig:Delaunay}
\end{figure}

We follow the construction in \cite{Aur90} to dualize the order-$k$
Voronoi tessellations.
For any $Q \subseteq X$, let $x_Q = \sum_{x \in Q} x / \card{Q}$ be the average point
with weight $w_Q = \norm{x_Q}^2 - \sum_{x \in Q} \norm{x}^2 / \card{Q}$.
The corresponding \emph{power function},
$\power{Q} \colon \Rspace^n \to \Rspace$, is defined by
$\power{Q} (p) = \Edist{p}{x_Q}^2 - w_Q$
and generalizes the squared Euclidean distance from $p$ to $x_Q$.
Let now $X_k$ be the collection of subsets $Q \subseteq X$ with $\card{Q} = k$.
The \emph{weighted Voronoi domain} of $Q \in X_k$ contains all points
$p \in \Rspace^n$ for which $\power{Q} (p) \leq \power{P} (p)$ for all $P \in X_k$,
and the (order-$1$) \emph{weighted Voronoi tessellation}
is the collection of non-empty such domains.
It can be proven that $p \in \domain{Q}$ iff $\power{Q} (p) \leq \power{P} (p)$
for all $P \in X_k$.
In other words, the order-$k$ Voronoi tessellation of $X$
is equal to the order-$1$ weighted Voronoi tessellation of $X_k$.
For the latter, there is a well-defined dual whose vertices are the
points $x_Q$ that have non-empty weighted domains.
It can be obtained as a projection of the lower convex hull of a special lifting of points
to $\Rspace^{n+1}$; see \cite{Aur87}.
We call this dual the \emph{order-$k$ Delaunay mosaic} of $X$, denoted $\Del{k}{X}$.
Figure \ref{fig:Delaunay} illustrates this construction by showing the
dual mosaics of the two Voronoi tessellations in Figure \ref{fig:Voronoi}.

We study $\Vor{k}{X}$ and $\Del{k}{X}$ when $X$ is a stationary
Poisson point process \cite{Kin93} with density $\density > 0$ in $\Rspace^n$.
With probability $1$, such a set $X$ is locally finite and in general position:
no $\DimDel+2$ points lie on a common $\DimDel$-plane
and no $\DimDel+3$ points line on a common $\DimDel$-sphere in $\Rspace^n$,
for $0 \leq \DimDel < n$.
The first result of this paper concerns the expected area of the
$\DimVor$-skeleton of an order-$k$ Poisson--Voronoi tessellation.
By definition, this is the $\DimVor$-dimensional Lebesgue measure of the
union of all $\DimVor$-dimensional faces of order-$k$ Voronoi domains.
Since this area is infinite, we normalize by letting $\Area{\DimVor}{k}{n}$
be the area of the $\DimVor$-skeleton within a unit volume of space.
\begin{theorem}[Expected Area]
  \label{thm:Area}
  Let $X$ be a stationary Poisson point process with density $\density > 0$
  in $\Rspace^n$, let $k\geq 1$ and $0 \leq \DimVor < n$.
  The expected area of the $\DimVor$-skeleton
  of the order-$k$ Voronoi tessellation of $X$ per unit volume of space is
  \begin{align}
    \EE{\Area{\DimVor}{k}{n}}  &=  \density^{\frac{n-\DimVor}{n}}
      \sum\displaylimits_{i=\max\{0, k+\DimVor-n\}}^{k-1}
        \frac{2^{n-\DimVor+1} \pi^{\frac{n-\DimVor}{2}}}
             {i! n (n - \DimVor + 1)!}
        \tfrac{\Gama{\frac{n^2-n\DimVor+\DimVor+1}{2}}
              \Gama{1+\frac{n}{2}}^{n-\DimVor+\frac{\DimVor}{n}}
              \Gama{n-\DimVor+i+\frac{\DimVor}{n}}}
             {\Gama{\frac{n^2-n\DimVor+\DimVor}{2}}
              \Gama{\frac{n+1}{2}}^{n-\DimVor}
              \Gama{\frac{\DimVor+1}{2}}} .
      \label{eqn:Area}
  \end{align}
\Skip{                        &=  \density^{\frac{n-\DimVor}{n}}
      \sum\displaylimits_{i=i_0}^{k-1}
        \frac{2^{n-\DimVor+1} \pi^{\frac{n-\DimVor}{2}}}
             {i! n (n - \DimVor + 1)!}
        \frac{\Gama{\frac{n^2-n\DimVor+\DimVor+1}{2}}
              \Gama{1+\frac{n}{2}}^{n-\DimVor+\frac{\DimVor}{n}}
              \Gama{n-\DimVor+i+\frac{\DimVor}{n}}}
             {\Gama{\frac{n^2-n\DimVor+\DimVor}{2}}
              \Gama{\frac{n+1}{2}}^{n-\DimVor}
              \Gama{\frac{\DimVor+1}{2}}} ,
      \label{eqn:Area}
  \end{align}
  in which $i_0 = \max \{ 0, k+\ell-n \}$.
}
  For $\DimVor = n$, we have $\EE{\Area{n}{k}{n}} = \Area{n}{k}{n} = 1$.
\end{theorem}

Our second result counts the cells in an order-$k$ Poisson--Delaunay mosaic.
Letting $G$ be a $\DimDel$-dimensional such cell,
we note that it uniquely determines the smallest sphere
centered at a point of the dual order-$k$ Voronoi polyhedron
such that at least $k$ points of $X$ lie inside or on the sphere;
see Section \ref{sec:4} for details.
We call the center and the radius of this sphere
the \emph{center} and the \emph{radius} of $G$.
To count, we specify a dimension $0 \leq \DimDel \leq n$,
a Borel region $\Region \subseteq \Rspace^n$,
and a radius $r_0 \geq 0$,
and we write $d_{\DimDel}^{k,n} (r_0)$ for the number
of $\DimDel$-cells in $\Del{k}{X}$ whose center
belongs to $\Region$ and whose radius is at most $r_0$.
We give an explicit formula for the expectation of
$d_{\DimDel}^{k,n} (r_0)$
using the constants $\CConst{p}{q}{n}$
defined in \cite{ENR16}.
\begin{theorem}[Expected Number of Cells]
  \label{thm:Number}
  Let $X$ be a stationary Poisson point process with density $\density > 0$
  in $\Rspace^n$, let $k \geq 1$ and $0 < \DimDel \leq n$.
  The expected number of $\DimDel$-cells in $\Del{k}{X}$ 
  with center in a Borel region $\Omega$ and radius at most $r_0$ satisfies
  \begin{align}
    \EE{d_{\DimDel}^{k,n}(r_0)}
      &= \density \norm{\Omega} \cdot \sum\displaylimits_{u=\DimDel}^{n}  
                   \sum\displaylimits_{v=1}^{u}
                     \CConst{v}{u}{n}
                   \sum\displaylimits_{g=1}^{g_1}
                     \frac{\iGama{u+k-g}{\density \nu_n r_0^n}}
                          {\Gama{k-g+1} \Gama{u}} 
                   \sum\displaylimits_{t=t_0}^{t_1}
                     \binom{v+1}{t} \binom{u-v}{t+\DimDel-v} ,
      \label{eqn:Number}
  \end{align}
  in which $g_1 = \min \{k,u\}$,
  $t_0 = \max\{0, v-\DimDel, g-\DimDel\}$, and
  $t_1 = \min\{v+1, u-\DimDel, g-1\}$.
	Further, for $\DimDel = 0$ and $k \geq 2$ we have
	\begin{align}
		\EE{\dNum{0}{k}{n} (r_0)}  &= 
			\density \norm{\Region} \cdot
			\sum\displaylimits_{u=1}^{n}  
			\sum\displaylimits_{v=1}^{u}
			  \CConst{v}{u}{n}
                          \frac{\iGama{u+k-v-1}{\density \nu_n r_0^n}}
                               {\Gama{k-v}\Gama{u}} .
	\end{align}	
\end{theorem}
Setting $r_0=\infty$, we obtain the expected total number of $\DimDel$-cells
in $\Del{k}{X}$ with center in $\Region$.
It is easy to verify that Theorem \ref{thm:Number} agrees
with \cite{ENR16} for $k=1$.
The case $\DimDel = 0$ is slightly different from the other dimensions;
and for $k = 1$ it is trivial because all points of $X$ are vertices of $\Del{1}{X}$.
Theorem \ref{thm:Number} implies that the radius of a \emph{typical $\DimDel$-cell} in $\Del{k}{X}$
follows a mixed Gamma distribution;
see \cite{ENR16}, where the details of this correspondence
are spelled out for the case $k = 1$.

Theorem	\ref{thm:Number} is derived as a corollary of the main technical
achievement of this paper:  the development of a discrete Morse theory
for order-$k$ Delaunay mosaics, and explicit formulas that count the
intervals in this theory.
Rather than presenting this result here,
we refer to Section \ref{sec:5} for its precise statement.

\ourparagraph{Outline.}
Section \ref{sec:2} describes the order-$k$ Voronoi tessellations
in detail, including a local characterization of their polyhedra
and a proof of Theorem \ref{thm:Area}.
Section \ref{sec:3} describes the order-$k$ Delaunay mosaics in detail,
including a complete classification of their cells.
Section \ref{sec:4} generalizes the discrete Morse theory of Delaunay mosaics
in \cite{BaEd17} from order-$1$ to order-$k$.
Section \ref{sec:5} counts the generalized intervals in 
the order-$k$ Delaunay mosaic, which leads to a proof of Theorem \ref{thm:Number}.
Section \ref{sec:6} concludes the paper.

\section{Voronoi Polyhedra}
\label{sec:2}

Any face of an order-$k$ Voronoi domain is a convex polyhedron that is
shared by a positive number of these domains.
Assuming its dimension is $\DimVor$, for some $0 \leq \DimVor \leq n$,
we call this face an \emph{order-$k$ Voronoi $\DimVor$-polyhedron}.
We begin with a geometric result about points on a sphere,
then use this result to prove a local characterization of
the order-$k$ Voronoi polyhedra, and finally prove Theorem \ref{thm:Area}.

\ourparagraph{Delaunay spheres.}
Let $X \subseteq \Rspace^n$ be locally finite.
For a point $p \in \Rspace^n$ and a positive integer $k$,
the \emph{order-$k$ Delaunay sphere} of $p$, denoted $\Sphere{k}{p}$,
is the smallest sphere centered at $p \in \Rspace^n$
such that the number of points of $X$ that lie inside or on the sphere is at least $k$.
To avoid possible ambiguities, we say a point lies \emph{inside} a sphere
if it belongs to the open ball bounded by the sphere.
It will be convenient to have short notation for these points
as well as their numbers.
Observing that $\interior{\conv{\Sphere{k}{p}}}$ is the open ball
with boundary $\Sphere{k}{p}$, we define
\begin{align}
  \Inn{k}{p} = X \cap \interior{\conv{\Sphere{k}{p}}}
    &\mbox{\rm ~~and~~} \hspace{2pt}\inn{k}{p} = \card{\Inn{k}{p}} ,                   \\
  \Onn{k}{p} = X \cap \Sphere{k}{p} ~~~~~~~~~~~
    &\mbox{\rm ~~and~~} \onn{k}{p} = \card{\Onn{k}{p}} .
\end{align}
By definition, $\inn{k}{p} + \onn{k}{p} \geq k$,
and by minimality of the radius, $\onn{k}{p} \geq 1$ and $\inn{k}{p} \leq k-1$.
The $\inn{k}{p}$ points in $\Inn{k}{p}$ are the unique $\inn{k}{p}$ nearest points to $p$,
the $\onn{k}{p}$ points in $\Onn{k}{p}$ are all at the same distance from $p$,
and all other points of $X$ are further from $p$.
With these notions, we get the following characterization of the order-$k$ Voronoi domains:
\begin{lemma}[Incident Voronoi Domains]
  \label{lem:local}
  Let $X \subseteq \Rspace^n$ be locally finite and in general position,
  and let $Q \subseteq X$ with $\card{Q} = k$.
  A point $p \in \Rspace^n$ belongs to $\domain{Q}$
  iff $\Inn{k}{p} \subseteq Q \subseteq \Inn{k}{p} \cup \Onn{k}{p}$.
\end{lemma}

\ourparagraph{Equivalence relation.}
We want to strengthen the previous lemma by
including polyhedra other than the Voronoi domains.
Recall that the interiors of the order-$k$ Voronoi polyhedra partition $\Rspace^n$.
To reconstruct this partition, we say that $p, q \in \Rspace^n$ are \emph{equivalent}
if their order-$k$ Delaunay spheres identify the same subsets of $X$.
More formally, we distinguish between the cases in which
the Delaunay sphere encloses $k$ or more than $k$ points:
\begin{align*}
  \! p \sim_X q  &\mbox{\rm ~if~}
    \left\{ \begin{array}{ll}
      \Inn{k}{p} \cup \Onn{k}{p} = \Inn{k}{q} \cup \Onn{k}{q} 
        &  \!\mbox{\rm for~} \inn{k}{p} + \onn{k}{p} = \inn{k}{q} + \onn{k}{q} = k , \\
      \Inn{k}{p} = \Inn{k}{q} , \Onn{k}{p} = \Onn{k}{q} 
        &  \!\mbox{\rm for~} \inn{k}{p} + \onn{k}{p} = \inn{k}{q} + \onn{k}{q} > k .
    \end{array} \right.
\end{align*}
We claim that the equivalence classes of $\sim_X$ are precisely the (relative) interiors
of the order-$k$ Voronoi polyhedra.
\begin{lemma}[Interiors of Order-$k$ Voronoi Polyhedra]
  \label{lem:equivalence}
  Let $X \subseteq \Rspace^n$ be locally finite and in general position.
  Then $p, q \in \interior{F}$, for an order-$k$ Voronoi polyhedron $F$,
  iff $p \sim_X q$.
\end{lemma}
\ourproof
  We first show that $p \sim_X q$ implies that the two points belong to the
  interior of a common order-$k$ Voronoi polyhedron.
  In the first case, when $\inn{k}{p} + \onn{k}{p} = \inn{k}{q} + \onn{k}{q} = k$,
  this is clear because $Q = \Inn{k}{p} \cup \Onn{k}{p} = \Inn{k}{q} \cup \Onn{k}{q}$
  is the unique set of $k$ nearest points in $X$,
  so $p, q \in \interior{\domain{Q}}$, which is an order-$k$ Voronoi $n$-polyhedron.
  In the second case, when $\inn{k}{p} + \onn{k}{p} = \inn{k}{q} + \onn{k}{q} > k$,
  we let $i = \inn{k}{p} = \inn{k}{q}$ and note that $i < k$.
  The points in $\Inn{k}{p} = \Inn{k}{q}$ are the unique $i$ nearest points,
  and we can add any $k-i$ points from $\Onn{k}{p} = \Onn{k}{q}$
  to get a complete set of $k$ nearest points.
  There are $\binom{\onn{k}{p}}{k-i} = \binom{\onn{k}{q}}{k-i}$ such choices,
  and by Lemma \ref{lem:local} each gives an order-$k$ Voronoi domain.
  These choices exhaust the domains that contain $p$ or $q$ on their boundaries.
  The set of points at equal distance from $\onn{k}{p} = \onn{k}{q}$ points
  of $X$ is a plane of dimension $n+1-\onn{k}{p} = n+1-\onn{k}{q}$,
  which implies that this is also the dimension of the order-$k$ Voronoi 
  polyhedron whose interior contains $p$ and $q$.

  We second show that $p \nsim_X q$ implies that $p$ and $q$ belong
  to the interiors of different order-$k$ Voronoi polyhedra.
  Assume the contrary. 
  We note that the dimension of the order-$k$
  Voronoi polyhedron whose interior contains $p$ is $n$,
  if $\inn{k}{p} + \onn{k}{p} = k$, and $n+1-\onn{k}{p}$,
  if $\inn{k}{p} + \onn{k}{p} > k$, and similar for $q$.
  In the first case, we would need
	$\inn{k}{q} + \onn{k}{q} = k$ to match the dimensions of the domains,
  but then $\Inn{k}{p} \cup \Onn{k}{p} \neq \Inn{k}{q} \cup \Onn{k}{q}$,
  so $p$ and $q$ belong to different domains.
  In the second case,
  we would need $\onn{k}{q} = \onn{k}{p}$ to have the same dimension of the polyhedra.
  Hence, $\Inn{k}{p} \neq \Inn{k}{q}$ or
  $\Inn{k}{p} = \Inn{k}{q}$ and $\Onn{k}{p} \neq \Onn{k}{q}$.
  In either case, we get a different collection of order-$k$ Voronoi domains
  for $p$ than for $q$.
\eop

\ourparagraph{Proof of Theorem \ref{thm:Area}.}
Recall that the proof of Lemma \ref{lem:equivalence} determines the dimension
of the order-$k$ Voronoi polyhedron whose interior contains a point $p \in \Rspace^n$ as
$n$, if $\inn{k}{p} + \onn{k}{p} = k$, and as
$n+1-\onn{k}{p}$, if $\inn{k}{p} + \onn{k}{p} > k$.
Equivalently, $p$ belongs to the interior of an order-$k$ Voronoi $\DimVor$-polyhedron iff
\begin{align}
  \DimVor = n              &\mbox{\rm ~and~} \inn{k}{p} + \onn{k}{p} = k \mbox{\rm ~~or}
    \label{eqn:first} \\
  0 \leq \DimVor \leq n-1  &\mbox{\rm ~and~} \onn{k}{p} = n-\DimVor+1 
                         \mbox{\rm ~and~} k+\DimVor-n \leq \inn{k}{p} \leq k-1.
    \label{eqn:second}
\end{align}
These relations suffice to extend the analysis in \cite{ScWe08}
from skeletons of order-$1$ to skeletons of order-$k$ Voronoi tessellations.
For $0 \leq \DimVor \leq n - 1$, they can be obtained as in
\cite[Theorem 10.2.4]{ScWe08}, which is the special case $k=1$ of \eqref{eqn:Area}.
The sole difference is that
we use the probability that there are $i$ points inside the sphere
instead of $0$, and sum over all admissible values of $i$, thus getting
$\Gama{n-\DimVor+i+\tfrac{\DimVor}{n}} / i!$ instead of
$\Gama{n-\DimVor+\tfrac{\DimVor}{n}}$ in the numerator.
This is precisely \eqref{eqn:Area}.
For $\DimVor = 0$ this gives the expected number of vertices in the order-$k$ Poisson--Voronoi mosaic.
For $\DimVor=n$ we trivially have $\EE{\Area{n}{k}{n}} = \Area{n}{k}{n} = 1$.
Theorem \ref{thm:Area} is thus proved.

\section{Delaunay Cells}
\label{sec:3}

In this section, we are more specific about the dual of the order-$k$ Voronoi tessellation.
As mentioned in Section \ref{sec:1}, each vertex of the order-$k$ Delaunay mosaic
is the average of the $k$ points that generate a non-empty order-$k$ Voronoi domain.
Each $(n-\DimDel)$-polyhedron of $\Vor{k}{X}$ is shared by a number of Voronoi domains,
each domain corresponds to a vertex, and the polyhedron corresponds to the $\DimDel$-cell
in $\Del{k}{X}$ that is the convex hull of these vertices.
Since $\Vor{k}{X}$ is not necessarily primitive, $\Del{k}{X}$ is not necessarily simplicial.

\begin{figure}[hbt]
  \centering \resizebox{!}{1.90in}{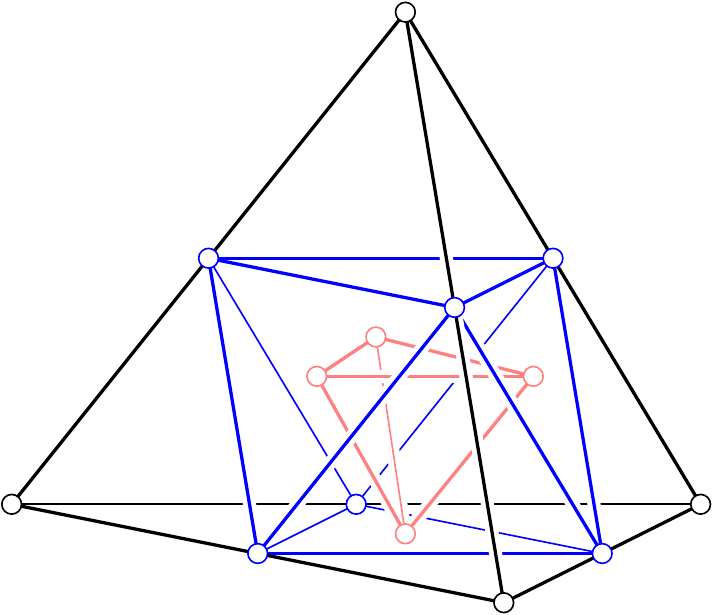}
  \caption{The three barycenter polytopes in $\Rspace^3$:
    the generation-$1$ tetrahedron, the generation-$2$ octahedron,
    and the generation-$3$ tetrahedron.}
  \label{fig:tetocta}
\end{figure}

\ourparagraph{Barycenter polytopes.}
We introduce a class of convex polytopes that is slightly richer than the class of simplices.
As we will see later, this class contains all polytopes we generically encounter in
order-$k$ Delaunay mosaics.
Let $\Delta^n$ be an $n$-dimensional simplex and recall that it has
$\binom{n+1}{g}$ faces of dimension $g-1$, for $1 \leq g \leq n + 1$.
The corresponding \emph{generation-$g$ barycenter polytope} is
the convex hull of the barycenters of all $(g-1)$-faces, denoted $\Delta_g^n$.
For $g = n+1$, the barycenter polytope is a single point,
but for other values of $g$ it is $n$-dimensional.
For $g = 1$ and $g = n$ the polytopes are $n$-simplices,
namely the convex hull of the $n+1$ vertices, $\Delta_1^n = \Delta^n$,
and the convex hull of the barycenters of the $n+1$ $(n-1)$-faces, $\Delta_n^n$.
For $2 \leq g \leq n-1$, the barycenter polytope is not a simplex,
and the first such case is $\Delta_2^3$, which is an octahedron;
see Figure \ref{fig:tetocta}.
A more detailed description of these polytopes is not needed,
and we refer to \cite{EdOs17} for additional information.

\ourparagraph{Characterization.}
If $X$ is in general position, which we assume,
then every cell of $\Del{k}{X}$ is a barycenter polytope.
To prove this, we consider a $u$-dimensional cell $G$ of $\Del{k}{X}$
and recall that all interior points
of its dual $(n-u)$-dimensional polyhedron $F$ of $\Vor{k}{X}$ are equivalent.
In other words, there are sets
$I = \Inn{k}{F}$ and $U = \Onn{k}{F}$ that uniquely determine $F$
as the polyhedron whose interior points $p$ satisfy $I = \Inn{k}{p}$ and $U = \Onn{k}{p}$.
We rewrite \eqref{eqn:first} and \eqref{eqn:second} to get constraints
on the sizes of the two sets:
\begin{align}
  \card{I} + \card{U} = k         
    &\mbox{\rm ~~~if~} u = 0,
    \label{eqn:first2}                                                \\
  \card{U} = u+1 \mbox{\rm ~and~} k-u \leq \card{I} \leq k-1
    &\mbox{\rm ~~~if~} u > 0.
    \label{eqn:second2}
\end{align}
The vertices of $\Del{k}{X}$ are governed by \eqref{eqn:first2},
while cells of higher dimensions are governed by \eqref{eqn:second2}.
Focusing on the cells of dimension $0 < u \leq n$, we note
that \eqref{eqn:second2} allows for a range of $u$ possible sizes of the set $I$.
These correspond to the generations of the barycenter polytopes,
as we now explain.
Let $i = \card{I}$ and define $g = k-i$,
noting that \eqref{eqn:second2} implies $1 \leq g \leq u$.
By Lemma \ref{lem:equivalence}, $F$ is the intersection of $\binom{u+1}{g}$ order-$k$ Voronoi domains
corresponding to $Q = I \cup \Uin$, in which $\Uin \subseteq U$ with $\card{\Uin} = g$.
So its dual cell $G$ is the convex hull of the averages $x_Q$ of these sets,
as discussed in Section \ref{sec:1}.
Writing each average as
\begin{align}
  x_Q  &=  \tfrac{1}{k} \left[ \sum\nolimits_{x \in I} x
                             + \sum\nolimits_{x \in \Uin} x \right]
        =  \tfrac{k-g}{k} x_I + \tfrac{g}{k} x_{\Uin} ,
  \label{eqn:avg}
\end{align}
we see that the convex hull of the $x_Q$ is a scaled and translated copy
of a generation-$g$ barycenter polytope,
namely the convex hull of the points $x_{\Uin}$.
Since $\card{U} = u+1$, this polytope is $u$-dimensional, as expected.
To summarize, we have a complete description of the cells in an order-$k$
Delaunay mosaic.
\begin{lemma}[Order-$k$ Delaunay Cells]
  \label{lem:Delf}
  Let $X \subseteq \Rspace^n$ be locally finite and in general position,
  and let $I, U \subseteq X$ with $I \cap U = \emptyset$.
  If $\card{I} + \card{U} = k$, then there is a point $p \in \Rspace^n$ with
  $\Inn{k}{p} \cup \Onn{k}{p} = I \cup U$ iff $x_{I \cup U}$ is a vertex of $\Del{k}{X}$.
  If $\card{I} + \card{U} \geq k+1$,
  then there is a point $p \in \Rspace^n$ with $\Inn{k}{p} = I$ and $\Onn{k}{p} = U$ iff
  the $u$-dimensional generation-$g$ barycenter polytope
  defined by $I$ and $U$ belongs to $\Del{k}{X}$,
  in which $u = \card{U} - 1$ and $g = k - \card{I}$.
\end{lemma}

\section{Relaxed Discrete Morse Theory}
\label{sec:4}

To count Delaunay cells in a stochastic setting, we would estimate
the probability that a given cell is defined by an order-$k$ Delaunay sphere.
For cells of intermediate dimension, there are pencils of possible
such spheres, which presents a challenge to the local methods of
probability theory.
To circumvent this difficulty, we follow the approach of \cite{ENR16}
and group the cells into intervals defined by a discrete Morse function;
see \cite{For98} for an introduction to discrete Morse theory,
and \cite{Fre09} for the generalization of the theory that fits
the geometry of Delaunay mosaics \cite{BaEd17}.
As we will see shortly,
order-$k$ Delaunay mosaics pose new difficulties,
which require a further relaxation of the theory.

\ourparagraph{Radius function.}
Recall that every $\DimDel$-cell $G \in \Del{k}{X}$ corresponds to
an $(n-\DimDel)$-polyhedron $F$ of $\Vor{k}{X}$.
By Lemma \ref{lem:equivalence}, 
for any point $p \in \interior{F}$, the Delaunay sphere $\Sphere{k}{p}$
passes through the same $\DimDel+1$ points $\Onn{k}{p} = \Onn{k}{F}$,
and $G$ is a scaled and translated
copy of a barycenter polytope defined by $\Onn{k}{p}$.
Since this is the smallest sphere centered at $p$
such that the number of points of $X$ that lie inside or on the sphere is at least $k$,
the sphere does not depend on $F$,
and its radius, $\radius{k}{p}$, is continuous as function of $p$.
Noting that $F$ is compact, we can therefore introduce
$\Rfun \colon \Del{k}{X} \to \Rspace$ defined by
\begin{align}
  \Rfun (G)  &=  \min \{ \radius{k}{p} \mid p \in F
                         \mbox{\rm ~and~} F \mbox{\rm ~dual to~} G \} ,
\end{align}
and call it the \emph{radius function} of $\Del{k}{X}$.
We call the point $p \in F$ that attains the minimum the \emph{center} of $G$.
This agrees with the definitions preceding Theorem \ref{thm:Number}.
Note that if the center $p$ of $G$ lies in the interior of a Voronoi face $F'$,
then $\Rfun (G')=\radius{k}{p}$ is the radius of $\Sphere{k}{p}$,
which determines the cell $G' \in \Del{k}{X}$ dual to $F'$
in the sense of Lemma \ref{lem:equivalence}.
An important observation is that $\Onn{k}{F} \subseteq \Onn{k}{p} = \Onn{k}{F'}$
and $\Inn{k}{F'} \subseteq \Inn{k}{F} \subseteq \Inn{k}{F'} \cup \Onn{k}{F'}$,
because all $k$-tuples of points of $X$, whose order-$k$ Voronoi domains intersect in $F$,
are involved in forming $F'$.
With this in mind, it is easy to determine which Voronoi polyhedra of any fixed dimension
contain $F'$.

\begin{figure}[hbt]
  \centering \resizebox{!}{1.7in}{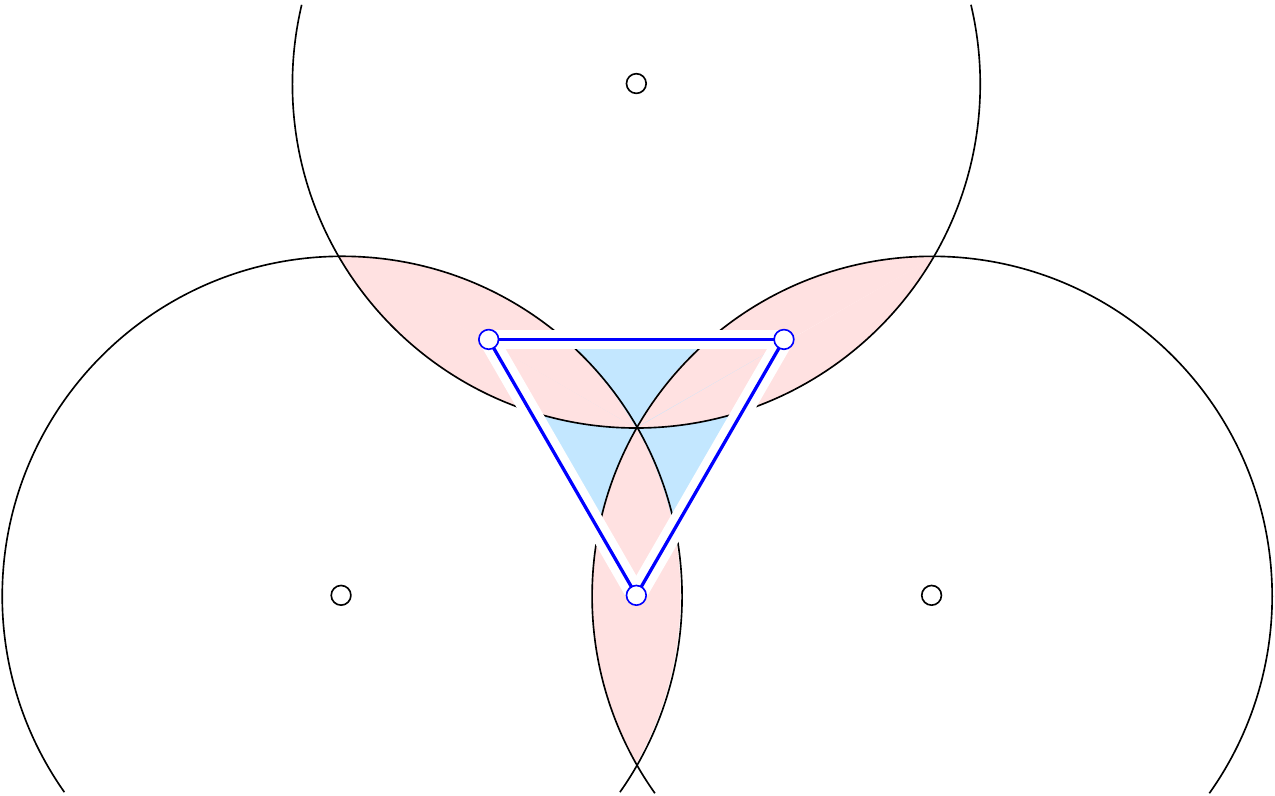}
  \caption{The radius function partitions the order-$2$ Delaunay mosaic
    of the three points into four relaxed intervals:
    three contain a vertex each, and the fourth relaxed interval
    contains the triangle together with its three edges.}
  \label{fig:counterexample}
\end{figure}
The discrete Morse theory of \cite{For98} requires that level sets of the radius function
are singletons and pairs,
while the generalized discrete Morse theory of \cite{Fre09}
allows intervals, which are maximal sets of faces
of a cell that share a common face.
The level sets of $\Rfun$ are not necessarily of this type,
as we now show.
Let $X$ consist of three points spanning an equilateral triangle
with unit length edges in the plane.
The order-$2$ Delaunay mosaic consists of
the triangle spanned by the midpoints of the three edges,
together with its edges and vertices.
Observe that $r_0 = 1/2$ is the radius assigned to its three vertices,
and $r_1 = \sqrt{3}/3$ is assigned to the triangle together with its
three edges; see Figure \ref{fig:counterexample}.
Indeed, the closed disks of radius $r$ centered at the points in $X$
have pairwise intersections iff $r \geq r_0$, and they have a
non-empty common intersection iff $r \geq r_1$.
Each vertex of $\Del{2}{X}$ has its own center in the
interior of the corresponding Voronoi $2$-polyhedron,
but the triangle and its three edges share the center
at the circumcenter of the triangle.
The triangle together with its edges is not an interval,
so $\Rfun$ is not a generalized discrete Morse function,
and we refer to it as a \emph{relaxed discrete Morse function}.
A justification of this terminology can be found at the end of this section.

\ourparagraph{Relaxed intervals.}
The radius function $\Rfun$ is \emph{monotonic}, by which we mean
that $\Rfun (G) \leq \Rfun (G')$ whenever $G$ is a face of $G'$.
However, equality is possible, namely when the order-$k$ Voronoi polyhedron $F'$
dual to $G'$
contains the center of $G$, which is in $F\setminus \interior{F}$.
By definition, a \emph{relaxed interval} of $\Rfun$ is a maximal
collection of cells in $\Del{k}{X}$ that share the center, and hence the function value.
Thus, every level set of $\Rfun$ is a disjoint union of relaxed intervals.

The previous example begs the question how much more general the relaxed intervals
are compared to the intervals.
Each relaxed interval has a unique \emph{upper bound}, which is a
cell $G \in \Del{k}{X}$, whose dual Voronoi polyhedron,
$F$, contains the center $p$ of $G$ in its interior.
Write $U = \Onn{k}{p}$ and $u = \card{U} - 1$.
The dimension of $G$ is thus $u$,
unless $\inn{k}{p}+\onn{k}{p} = k$, in which case it is $0$.
Considering any partition of $U$ into three sets, $U = \Uin \cup \Uon \cup \Uout$
with $\Uon \ne U$,
we can slightly perturb the sphere $\Sphere{k}{p}$ into a sphere $\SphereOnly$
such that $\Uon$ and $\Inn{k}{p} \cup \Uin$ are the points on and inside $\SphereOnly$,
respectively.
If the sizes of these two sets satisfy the requirements for the order-$k$ Delaunay sphere, 
they define a cell of $\Del{k}{X}$, which is a face of $G$.
On the other hand, every face of $G$ \emph{induces} such a partition.
We therefore get a correspondence between such partitions of $U$
and the faces of $G$, which is one-to-one unless $\card{\Uin \cup \Uon \cup \Inn{k}{p}} = k$,
in which case we get the same vertex for all partitions with the same $\Uout$.

We are particularly interested in distinguishing
the faces that share the center, $p$, from the other faces of $G$.
The crucial concept is the visibility of facets of $\conv{U}$ from $p$, which we now introduce.
Recall that the center $p$ of $G$ lies in the interior
of the dual Voronoi polyhedron by assumption.
It follows that $p$ lies in the affine hull of $U$.
Equivalently, the $u$-sphere with center $p$ that passes through
the $u+1$ points of $U = \Onn{k}{p}$ is a great-sphere of $\Sphere{k}{p}$.
The convex hull of $U$ is a $u$-simplex with
$u+1$ $(u-1)$-dimensional faces, which we call its \emph{facets}.
A facet is \emph{visible} from $p$ if the affine hull of the facet,
which is a $(u-1)$-plane, separates $p$ from $\conv{U}$
within the affine hull of $U$, which is a $u$-plane.
Let $v$ be the number of invisible facets minus $1$ and observe that
$v \geq 1$ because the $u+1$ points of $U$ lie on a sphere around $p$.
Let $V \subseteq U$ contain the points that belong to all visible facets,
and observe that $\card{V} = v+1$ because a vertex belongs to $V$
iff the facet opposite to the vertex is invisible. In particular, $V=U$ if there are no visible facets.
With these notions, we can identify the partitions of $U$
that correspond to faces of $G$ in the relaxed interval with upper bound $G$.
\begin{lemma}[Visibility and Relaxed Intervals]
  \label{lem:vis_int}
  Let $X \subseteq \Rspace^n$ be locally finite and in general position.
  Let $G \in \Del{k}{X}$ with corresponding order-$k$ Delaunay sphere $\Sphere{k}{p}$
  be the upper bound of a relaxed interval of the radius function.
  A face $G'$ of $G$ belongs to the same relaxed interval
  iff the partition $\Onn{k}{p} = \Uin \cup \Uon \cup \Uout$
  induced by $G'$ satisfies $\Uin \subseteq V \subseteq \Uin \cup \Uon$.
\end{lemma}
\ourproof
  Write $U = \Onn{k}{p}$.
  Let $q$ be the center of $G'$,
  and recall that $G, G'$ belong to the same relaxed interval iff $p = q$.
  We have $p \neq q$ unless the following two conditions hold:
  \begin{description}
    \item[{\rm (i)}]  If an invisible face of $\conv{U}$ contains $\Uon$,
      then the opposite vertex must be in $\Uin$.
    \item[{\rm (ii)}]  If a visible face of $\conv{U}$ contains $\Uon$,
      then the opposite vertex must be in $\Uout$.
  \end{description}
  To see (i), we would move the center, $p$, normal to and slightly toward the
  facet while adjusting the radius so the sphere keeps passing through all vertices
  of the facet.
  This generates a smaller sphere for the same partition of $U$, hence $p \neq q$.
  The symmetric argument proves (ii).
  Now (i) is equivalent to $\Uin \subseteq V$,
  and (ii) is equivalent to $U \setminus V \subseteq \Uout$.
  Hence $p = q$ implies $\Uin \subseteq V \subseteq \Uin \cup \Uon$.
  The converse is also true because the two conditions prohibit a smaller
  sphere in the normal directions of all facets.
  These directions span all directions in the affine hull of $U$.
\eop
The only case when the induced decomposition is not necessarily unique,
is when $G'$ is a vertex.
In particular, if the upper bound $G$ is a vertex itself,
then we get $V=U$ as an additional requirement.

\ourparagraph{Critical and non-critical cases.}
We call a case \emph{critical} if the defining simplex, $\conv{U}$,
has no visible facets, and we call it \emph{non-critical} otherwise.
This classification is motivated by constructing $\Del{k}{X}$ incrementally,
adding one relaxed interval at a time in the order of the radius function.
In the critical case, the effect of adding the cells in the relaxed interval
changes the homotopy type of the current complex,
while in the non-critical case the homotopy type remains unchanged.
The proof of this claim is beyond the scope of this paper
and can be found in the yet unpublished \cite{EdOs17}.
Indeed, we are primarily interested in the number of cells per relaxed interval,
but the mentioned topological fact justifies that we call $\Rfun$
a relaxed discrete Morse function and not something much more general.

\section{Counting}
\label{sec:5}  

In this section, we count the cells in the relaxed intervals that arise
in the partition of order-$k$ Delaunay mosaics.
We then use the result to prove Theorem \ref{thm:Number}.

\ourparagraph{Cells in relaxed intervals.}
As explained in Section \ref{sec:4}, every relaxed interval has a unique upper bound,
which is a cell $G \in \Del{k}{X}$ whose center, $p \in \Rspace^n$,
is contained in the interior of the dual Voronoi polyhedron.
The order-$k$ Delaunay sphere of this point, $\Sphere{k}{p}$,
completely determines $G$; see \eqref{eqn:avg}.
Ignoring the case in which $G$ is a vertex,
we assume that $\inn{k}{p} + \onn{k}{p} \geq k+1$,
in which case $u = \onn{k}{p} - 1 \geq 1$ is the dimension of $G$
and $g = k - \inn{k}{p}$ is its generation.
As discussed above, different vertices of $G$  correspond to
different subsets $\Uout$ of $U = \Onn{k}{p}$ with $\card{\Uout} = \card{U} - g$.
To get the number of vertices, we therefore count
the partitions $U = \Uin \cup \Uout$ with $\card{\Uin} = g$; compare with \eqref{eqn:first2}.
To get the number of $\DimDel$-faces of $G$ for $0 < \DimDel \leq u$,
we count the partitions $U = \Uin \cup \Uon \cup \Uout$ that satisfy
$\card{\Uon} = \DimDel+1$ and $g-\DimDel \leq \card{\Uin} \leq g-1$; compare with \eqref{eqn:second2}.
To further limit the number to the cells in the relaxed interval of $G$,
we use Lemma \ref{lem:vis_int} and restrict to $\Uin \subseteq V \subseteq \Uin \cup \Uon$,
in which $V \subseteq U$ with $\card{V} = v + 1$ contains the vertices that belong to all visible
facets of $U$.

For $\DimDel = 0$, the last condition is equivalent to $\Uin = V$.
Writing $\NFaces{v}{g}{u}{j}$ for the number of faces in the relaxed interval with
upper bound $G$, we therefore have $\NFaces{v}{g}{u}{0} = 1$ if $g = v + 1$,
and $\NFaces{v}{g}{u}{0} = 0$ otherwise.
When $\DimDel > 0$, the dimension requirement is that $\card{\Uon} = \DimDel + 1$.
Writing $t = \card{\Uin}$, we can formulate the question purely combinatorially,
first choosing the union $\Uin \cup \Uon \subseteq U$ such that $V \subseteq \Uin \cup \Uon$ and
second choosing $\Uin \subseteq V$:
\emph{how many ways are there to pick $(t + \DimDel + 1) - (v + 1)$ from $(u+1) - (v+1)$ points
      and then $t$ from $v+1$ points?}
The answer gives the number of faces in the relaxed interval:
\begin{align}
  \NFaces{v}{g}{u}{\DimDel}  &=  \sum_{t = t_0}^{t_1} \binom{u - v}{t + \DimDel - v} \binom{v + 1}{t},
  \label{eqn:NFaces}
\end{align}
in which $t_0 = \max \{0, v-\DimDel, g-\DimDel \}$ and $t_1 = \min \{ v+1, u-\DimDel, g-1 \}$
are obtained from $0 \leq t \leq v+1$, $0 \leq \DimDel-v+t \leq u-v$,
and $g-\DimDel \leq t \leq g-1$.
The first two conditions assert that the binomial coefficients make sense,
while the last one is the geometric requirement for the number of points inside the sphere.

\ourparagraph{Determination of intervals.}
The analysis in the previous section suggests
we use the order-$k$ Delaunay spheres as intrinsic characterization of the relaxed intervals.
Let $U \subseteq X \subseteq \Rspace^n$ with $\card{U} = u + 1 \leq n + 1$ be a simplex,
such that there are between $k-u-1$ and $k-1$ points inside
the smallest circumscribed sphere $\SphereOnly$ of $U$.
Letting $p$ be the center of this sphere,
we notice that $\SphereOnly = \Sphere{k}{p}$ and $\Onn{k}{p} = U$.
If $\onn{k}{p} + \inn{k}{p} > k$, it defines a cell $G$ of $\Del{k}{X}$,
namely a barycenter polytope of type $\Delta_g^u$, for $g = k - \inn{k}{p}$.
By Lemma \ref{lem:vis_int}, this cell is the upper bound of a relaxed interval 
of the radius function $\Rfun$,
which contains all cells that share $p$ as their center.
The lemma also asserts that the interval is
fully described by the set of vertices of $U$ that belong to all visible facets.
Writing $V$ for this set and $v = \card{V} - 1$ for its dimension,
we call $(v, u, g)$ the \emph{type} of the relaxed interval.
It is fully defined by $U$.

If $\onn{k}{p} + \inn{k}{p} = k$, then $p$ belongs
to the interior of the order-$k$ Voronoi domain of $\Onn{k}{p} \cup \Inn{k}{p}$.
By Lemma \ref{lem:vis_int} and the remark after it,
$p$ is the center of this domain iff it lies in the interior of $\conv{U}$.
In this case, we get a critical vertex, with $V = U$ and $g = u + 1$.
The type of this interval is thus $(u, u, u + 1)$.
This should not be confusing because vertices with different
relaxed interval types are really different kinds of vertices in the mosaic.

\ourparagraph{Proof of Theorem \ref{thm:Number}.}
We now apply the developed theory to prove our second main result.
Let $X$ be a stationary Poisson point process with density $\density > 0$ in $\Rspace^n$.
Using the intrinsic characterization,
we want to compute the expected numbers of intervals of type $(v, u, g)$,
while restricting the radius from above.
Write $\cNumNB{v}{u}{g}{k}{n} (r_0)$ for the number of 
tuples of $u+1$ points in $X$, whose
smallest circumspheres have $k-g$ points inside, have their center
in some region $\Region \subseteq \Rspace^n$, and have radius at most $r_0$.
As the previous discussion shows,
it is the same as the number $\cNum{v}{u}{g}{k}{n} (r_0)$
of intervals of type $(v, u, g)$ with center in $\Region$ and radius at most $r_0$,
when $1 \leq g \leq \max\{u, k\}$ or $v = u = g - 1$.
Following the approach in \cite{ENR16}, we focus on the non-trivial case $u > 0$ and use the 
Slivnyak--Mecke formula to express the expectation of this number as
\begin{align}
  \EE{\cNumNB{v}{u}{g}{k}{n} (r_0)}  &=
    \tfrac{1}{(u+1)!} \int\displaylimits_{\xxx \in (\Rspace^n)^{u+1}}
      \One{\Region} (\xxx) \One{r_0} (\xxx) \One{u-v} (\xxx)
      \PP{k-g}{\xxx} \density^{u+1} \diff \xxx,
\end{align}  
in which $\nu_n$ is the volume of a unit $n$-ball,
$\PP{k-g}{\xxx} = (\density \nu_n r^n)^{k-g} e^{-\density \nu_n r^n} / (k-g)!$
is the probability that the smallest circumsphere of $\xxx$ has $k-g$ points of $X$ inside,
$\One{\Region} (\xxx)$ indicates whether the center of this sphere
belongs to $\Region$,
$\One{r_0} (\xxx)$ indicates whether its radius is at most $r_0$,
and $\One{u-v} (\xxx)$ indicates whether $\xxx$ has $u-v$ visible facets.
The notation we use mimics the one in \cite{ENR16}, in particular,
we write $\xxx$ for a sequence of $u+1$ points,
which is better suited for integration than the set $U$ of $u+1$ points.
The only difference to Equation (3.4) in \cite{ENR16} is the use of
$\PP{k-g}{\xxx}$ instead of $\PP{0}{\xxx} = \PP{\emptyset}{\xxx}$.
As explained in that article, we can use the spherical Blaschke--Petkantschin
formula to compute this integral,
and to avoid redundancy, we focus on the differences.
Specifically, instead of
$\int\displaylimits_{r = 0}^{r_0} r^{nk-1} e^{-\density r^n \nu_n}
  = \tfrac{\iGama{k}{\density \nu_n r_0^n}}{n(\density \nu_n)^k}$
in (3.6) of \cite{ENR16}, we have
\begin{align*}
  \int\displaylimits_{r = 0}^{r_0} r^{n u-1}
                                      \frac{(\density r^n \nu_n)^{k-g}}{(k-g)!}
                                      e^{-\density r^n \nu_n}
    &=  \frac{(\density \nu_n)^{k-g}}{(k-g)!}
        \frac{\iGama{u+k-g}{\density \nu_n r_0^n}}{n(\density \nu_n)^{u+k-g}}
     =  \frac{\iGama{u+k-g}{\density \nu_n r_0^n}}{(k-g)! n (\density \nu_n)^{u}}.    
\end{align*}
Arguing exactly like in Lemma 3.1 in \cite{ENR16}, we get
\begin{align}
  \EE{\cNumNB{v}{u}{g}{k}{n} (r_0)}
    &=  \frac{\iGama{u+k-g}{\density \nu_n r_0^n}}{(k-g)!\Gama{u}}
        \CConst{v}{u}{n} \cdot \density \norm{\Region},
    \label{eqn:types}
\end{align}
in which the constant $\CConst{v}{u}{n}$ is as defined in \cite{ENR16}.
The case $u=0$ is exceptional, because the smallest circumscribed sphere of any single vertex
has radius $0$ and no points inside, so the only non-zero value is
$\EE{\cNumNB{0}{0}{1}{1}{n} (r_0)} = \density \norm{\Region}$
for all $r_0 \geq 0$, independent of the radius.
Returning to the number of relaxed intervals, we thus have
$\EE{\cNum{v}{u}{g}{k}{n} (r_0)} = \EE{\cNumNB{v}{u}{g}{k}{n} (r_0)}$
for admissible values of parameters,
i.e., for $1 \leq g \leq \min \{k, u\}$ or $v = u = g - 1$, and $0$ otherwise.
The result agrees with \cite{ENR16} for $k=1$.

Now that we have expressions for the number of relaxed intervals
of all types, it is not difficult to count the $\DimDel$-cells
in the order-$k$ Delaunay mosaic whose value under the radius function
is at most $r_0$:
\begin{align}
  \EE{\dNum{\DimDel}{k}{n} (r_0)}  &= 
    \sum\displaylimits_{u=\DimDel}^{n}
    \sum\displaylimits_{v=0}^{u}
    \sum\displaylimits_{g=1}^{\min\{k, u + 1\}}
      \NFaces{v}{g}{u}{\DimDel} \cdot
      \EE{\cNum{v}{u}{g}{k}{n} (r_0)} .
\end{align}
For $\DimDel>0$, we  can use \eqref{eqn:NFaces} and \eqref{eqn:types} to get
\begin{align}
  \EE{\dNum{\DimDel}{k}{n} (r_0)}  &= 
    \sum\displaylimits_{u=\DimDel}^{n}  
    \sum\displaylimits_{v=1}^{u}
    \sum\displaylimits_{g=1}^{g_1}
    \sum\displaylimits_{t=t_0}^{t_1}
      \binom{v+1}{t} \binom{u-v}{t+\DimDel-v}
      \frac{\iGama{u+k-g}{\density \nu_n r_0^n}}
           {(k-g)! \Gama{u}}
      \CConst{v}{u}{n}
      \cdot \density \norm{\Region} ,
\end{align}
in which $g_1 = \min\{k,u\}$, $t_0 = \max \{ 0, v-\DimDel, g-\DimDel\}$, and
$t_1 = \min \{ v+1, u-\DimDel, g-1\}$, as before.
For $\DimDel=0$ and $k \geq 2$, we take the sum of the numbers of relaxed intervals
with $g = v + 1$:
\begin{align}
  \EE{\dNum{0}{k}{n} (r_0)}  &= 
    \sum\displaylimits_{u=1}^{n}  
    \sum\displaylimits_{v=1}^{u}
      \frac{\iGama{u+k-v-1}{\density \nu_n r_0^n}}
           {(k-v-1) !\Gama{u}}
    \CConst{v}{u}{n}
    \cdot \density \norm{\Region}.
  \label{eq:null}
\end{align}
This completes the proof of Theorem \ref{thm:Number}.

\section{Discussion}
\label{sec:6}  

This paper gives evidence of the power of the discrete Morse theory
approach to questions in stochastic geometry.
The first step is the relaxation of discrete Morse functions
so they apply to order-$k$ Delaunay mosaics.
This relaxation is non-trivial and of independent interest.
Here we provide a complete combinatorial analysis of the relaxed intervals
that make up the discrete theory,
and we use it to generalize the main stochastic relations
in \cite{ENR16} from order-$1$ to order-$k$ Delaunay mosaics.

While the results in this paper are predominantly combinatorial
and probabilistic, there are connections to other areas of mathematics
and to applications outside of mathematics.
Results about the topological meaning of the relaxed
Morse theory are under investigation in the forthcoming \cite{EdOs17},
including algorithms to compute the persistent homology of multi-covers
with balls.
We hope that the stochastic and the topological tools together
give a novel approach to dealing with dense data
and will lead to a refined understanding of medium- to long-range
effects in locally finite configurations,
as they arise for example during the emergence
of order in particle arrangements.


\end{document}